\documentclass[amscd,amssymb,verbatim,a4paper]{amsart}

\theoremstyle{plain}
\newtheorem{theorem}{Theorem}[section]
\newtheorem{lemma}[theorem]{Lemma}
\newtheorem{proposition}[theorem]{Proposition}
\newtheorem{corollary}[theorem]{Corollary}

\theoremstyle{definition}

\theoremstyle{remark}
\newtheorem{remark}[theorem]{Remark}

\newcommand{\ic}{\ensuremath{\mathcal{I}}}
\newcommand{\oc}{\ensuremath{\mathcal{O}}}

\newcommand{\fc}{\ensuremath{\mathcal{F}}}
\newcommand{\dc}{\ensuremath{\mathcal{D}}}

\newcommand{\hc}{\ensuremath{\mathcal{H}}}
\newcommand{\cc}{\ensuremath{\mathcal{C}}}
\newcommand{\lc}{\ensuremath{\mathcal{L}}}

\newcommand{\Sc}{\ensuremath{\mathcal{S}}}
\newcommand{\tc}{\ensuremath{\mathcal{T}}}
\newcommand{\nc}{\ensuremath{\mathcal{N}}}

\newcommand{\pc}{\ensuremath{\mathcal{P}}}

\newcommand{\mc}{\ensuremath{\mathcal{M}}}
\newcommand{\qc}{\ensuremath{\mathcal{Q}}}
\newcommand{\rc}{\ensuremath{\mathcal{R}}}

\newcommand{\Pcq}{\mathbb{P}^5}

\newcommand{\Ps}{\mathbb{P}}

\newcommand{\bC}{\mathbb{C}}

\def\bin #1#2 {\left( \matrix { #1 \cr #2 \cr } \right) }

\newcommand{\tG}{\tilde{G}}

\newcommand{\tp}{\tilde{\pc}}

\newcommand{\dg}{\tilde{D}(g)}
\newcommand{\df}{\tilde{D}(f)}
\newcommand{\dtg}{\tilde{\mathbb{D}}(g)}
\newcommand{\dtf}{\tilde{\mathbb{D}}(f)}
\newcommand{\ddg}{{\mathbb{D}}(g)}
\newcommand{\ddf}{{\mathbb{D}}(f)}

\begin{document}

\title[Factoriality and N\'eron-Severi  groups]
{Factoriality and N\'eron-Severi  groups of a projective
codimension two complete intersection with isolated singularities}

\author{Vincenzo Di Gennaro }
\address{Universit\`a di Roma \lq\lq Tor Vergata\rq\rq, Dipartimento di Matematica,
Via della Ricerca Scientifica, 00133 Roma, Italy.}
\email{digennar@axp.mat.uniroma2.it}

\author{Davide Franco }
\address{Universit\`a di Napoli
\lq\lq Federico II\rq\rq, Dipartimento di Matematica e
Applicazioni \lq\lq R. Caccioppoli\rq\rq, P.le Tecchio 80, 80125
Napoli, Italy.} \email{davide.franco@unina.it}

\date{12/05/2006}

\bigskip
\abstract For a projective variety $Z$ and for any integer $p$,
define the $p$-th N\'eron-Severi  group $NS_p(Z)$ of $Z$ as the
image of the cycle map $A_{p}(Z)\to H_{2p}(Z; \mathbb{C})$. Now
let $X\subset \Ps^{2m+1}$ ($m\geq 1$) be a  projective variety of
dimension $2m-1$, with isolated singularities, complete
intersection of a smooth hypersurface of degree $k$, with a
hypersurface of degree $n>max\{k, 2m+1\}$, and let $F$ be a
general hypersurface of degree $n$ containing $X$. We prove that
the natural map $NS_m(X)\to NS_m(F)$ is surjective, and  that if
$dim\,NS_m(F)=1$ then $dim\,NS_m(X)=1$. In particular
$dim\,NS_m(X)=1$ if and only if $dim\,NS_m(F)=1$. When $X$ is a
threefold (i.e. $m=2$) we deduce a new characterization for the
factoriality of $X$, i.e. that $X$ is factorial if and only if
$dim\, NS_2(F)=1$. This allows us to give  examples of factorial
threefolds, in some case with many singularities. During the proof
of the announced results, we show that the quotient of the middle
cohomology of $F$ by the cycle classes coming from $X$ is
irreducible under the monodromy action induced by the
hypersurfaces of degree $n$ containing $X$. As consequences we
deduce a Noether-Lefschetz Theorem for a projective complete
intersection with isolated singularities, and, also using a recent
result on codimension two Hodge conjecture, in the case $X\subset
\Ps^{5}$ is a threefold as before,  we deduce that the general
hypersurface $F$ of degree $n$ containing $X$ verifies Hodge
conjecture.



\endabstract

\maketitle

\section{Introduction}
Let $X\subset \Pcq$ be a complex, projective, complete
intersection threefold, with isolated singularities. One says that
$X$ is {\it{factorial}} if its graded ring is unique factorization
domain. This is equivalent to the fact that every surface lying in
$X$ is a complete intersection on $X$ (\cite{Greco}, p. 69). Using
Lefschetz type Theorems (\cite{BS}, pg. 50-51), one sees that
factoriality is  also equivalent to say that every Weil divisor of
$X$ is a Cartier divisor, that it  is also equivalent to say that
$X$ is {\it{$\mathbb{Q}$-factorial}}, i.e. that every Weil divisor
of $X$ has a multiple which is a Cartier divisor (\cite{BS}, pg.
5-6), and that it is also equivalent to say that $X$ is {\it
locally factorial} (\cite{Greco}, pg. 69). For the interest in
this notion in the study of birational geometry, we refer to
\cite{KM}, \cite{MP}, \cite{Cheltsov1} and \cite{Ml}.

If $X$ is nonsingular then it is  factorial. This is no longer
true when $X$ is singular. For instance, when $X$ is a complete
intersection of general hypersurfaces containing a fixed plane,
then $Sing(X)$ consists of $(n+k-2)^2-(n-1)(k-1)$ ordinary double
points, where $(k,n)$ is the bi-degree of $X$, and $X$ is not
factorial for it contains the given plane. However it is known
that if $X$ is a complete intersection on a {\it smooth} fourfold
$G\subset\Pcq$ of degree $k$ with a hypersurface of degree $n\geq
k$ (e.g. $X\subset \mathbb{P}^4$), and if $X$ presents few
ordinary double points with respect to $k$ and $n$, then it is
factorial (see \cite{Cheltsov1}, \cite{Cheltsov2},
\cite{Cheltsov3}, \cite{Cheltsov4}, \cite{CD1}, \cite{CD2}).
Actually one conjectures that previous number is the sharp bound,
i.e. that any $X$ as above, with only ordinary double points and
such that $|Sing(X)|<(n+k-2)^2-(n-1)(k-1)$, is factorial.

When all the singularities of $X$ are   ordinary double points,
i.e. are   nodes, an important invariant in the study of the
birational geometry of $X$, introduced by C.H. Clemens (see
\cite{Clemens}, \cite{Werner}), is the defect $\delta(X)$ of $X$,
which is equal to the difference between the fourth and the second
Betti number of $X$. If $X$ is a complete intersection, the second
Betti number is equal to $1$ and so we have
$$
rk\, H_4(X; \mathbb{Z})=1+\delta(X).
$$
If in addition $X$ is a complete intersection on a smooth fourfold
$G\subset \Pcq$ of degree $k$ with a hypersurface of degree $n$
(e.g. $X\subset \mathbb{P}^4$), then by \cite{Cynk} we know that
$\delta(X)$ is the number of dependent conditions that vanishing
at singularities of $X$ imposes on the global sections of the line
bundle $\mathcal{O}_{G}(2n+k-6)$ on $G$. Using the same argument
as in \cite{Cheltsov2}, Remark 19, one deduces that $X$ is
factorial if and only if its defect vanishes, i.e. in the nodal
case one has:
\begin{equation}
\label{fact} {\text{$X$ is factorial if and only if $rk\, H_4(X;
\mathbb{Z})=1,$}}
\end{equation}
(see \cite{Cheltsov2}, \cite{Cheltsov3}, \cite{Cheltsov4},
\cite{Ml}). In particular we see that the factoriality can be
described  by a global topological condition, and that it depends
on the position of the nodes in the projective space. Without
assuming that the isolated singularities are nodes,
characterization (\ref{fact}) holds true only in one direction,
i.e.
\begin{equation}
\label{facto} {\text{if $rk\, H_4(X; \mathbb{Z})=1$ then $X$ is
factorial.}}
\end{equation}
In fact, if $rk\, H_4(X; \mathbb{Z})=1$ then  any projective
integral surface $S$ contained in $X$ has a general hyperplane
section $S\cap H$ which is homologous to a multiple of the
hyperplane section of $X\cap H$. Hence by Hamm-Lefschetz Theorem
it follows that $S\cap H$ is a complete intersection on $X\cap H$,
and so  $S$ is on $X$. As we said before, this means that $X$ is
factorial. On the other hand, from Noether-Lefschetz Theorem, the
cone $X\subset\Ps^4$ over a general surface $S\subset \Ps^3$ of
degree $n\geq 4$ is factorial, but $H_4(X; \mathbb{Z})\simeq
H_2(S; \mathbb{Z})$ (see \cite{Dimca}, p. 169, (4.18)) and the
last group has rank $> 1$. Therefore the converse of property
(\ref{facto}) is false. Notice that property (\ref{facto}) allows
us to give  examples of factorial non-nodal threefolds $X$. In
fact, by \cite{Dimca}, Theorem (4.17), we know that {\it{if
$X\subset \mathbb{P}^4$ has just one ordinary singular point of
multiplicity $<deg(X)$, then $rk(H_4(X; \mathbb{Z}))=1$}}.

Using the notion of N\'eron-Severi group, one may reformulate
property (\ref{facto}) in the following way. First recall that for
a projective variety $Z$ and for any integer $p$, one may define
the $p$-th N\'eron-Severi  group $NS_p(Z)$ of $Z$ as the image of
the cycle map $A_{p}(Z)\to H_{2p}(Z; \mathbb{C})$ (see
\cite{Fulton}, Chapter 19). Next recall that for a projective
complete intersection $Z$ of dimension $2m+\epsilon-1$, $0\leq
\epsilon \leq 1$, with isolated singularities, the only
interesting N\'eron-Severi group is $NS_m(Z)$ (see \cite{Dimca},
p. 161, Theorem (4.3)). Now notice that previous argument in
proving property (\ref{facto}) works well also if one simply
assumes that $dim(NS_2(X))=1$, and so, for a threefold $X$
complete intersection with isolated singularities, we have:
\begin{equation}
\label{factor} {\text{$X$ is factorial if and only if
$dim(NS_2(X))=1$. }}
\end{equation}

In the present paper we show   a new characterization of the
factoriality. Roughly saying, we prove that {\it $X$ is factorial
if and only if it is a complete intersection on a smooth fourfold
$F$ such that $dim(NS_2(F))=1$}. More precisely, first we prove
the following general result:

\medskip
\begin{theorem}
\label{main} Let $X\subset \Ps^{2m+1}$  be a projective variety of
dimension $2m-1\geq 1$, complete intersection of two smooth
hypersurfaces $F$ and $G$ of degrees $n$ and $k$, with $n>k$. Then
$X$ has isolated singularities. Moreover, if  $dim\, NS_m(F)=1$
 then $dim\, NS_{m}(X)=1$.
\end{theorem}
\medskip

The property that a complete intersection of two smooth
hypersurfaces of different degrees has at most isolated
singularities is a direct consequence of Proposition 4.3.6. in
\cite{Vogel}, and it holds true also when intersecting two smooth
hypersurfaces in a projective space of {\it even} dimension (see
also \cite{Lazarsfeld}, Example 6.3.8.). On the contrary, in the
projective space $\Ps^{2m+2}$ the assertion {\it if $dim\,
NS_m(F)=1$ then $dim\, NS_{m}(X)=1$} is false. In fact, for any
smooth hypersurface $F$ of odd dimension, one has $dim\,
NS_m(F)=1$, but there exist smooth complete intersections $X$ of
dimension $2m$ with $dim\, NS_{m}(X)>1$.

Even if one assumes that $X$ has isolated singularities, the
hypothesis $G$ smooth in Theorem \ref{main} is necessary, as
Example 5 in \cite{Cheltsov2} proves. Next we prove that,
conversely,  the hypotheses that $F$ is smooth and that $dim\,
NS_m(F)=1$ are also necessary for the property $dim\,
NS_{m}(X)=1$, in the following sense:

\medskip
\begin{theorem}
\label{mmain} Let $X\subset \Ps^{2m+1}$  be a projective variety
of dimension $2m-1\geq 1$, with isolated singularities, complete
intersection of a smooth hypersurface $G$ of degree $k$, with a
hypersurface of degree $n$. Set $\lc =\mid \ic _{X,\Ps^{2m+1}}(n)
\mid $, and let $F\in\lc$ be a general hypersurface. If $n\geq k$
then $F$ is smooth. Moreover, if $n>max\{k, 2m+1\}$  then the
natural map $NS_m(X) \to NS_m(F)$ is surjective.
\end{theorem}
\medskip

We will see that the analogous claim of Theorem \ref{mmain} in a
projective space of {even} dimension remains true, but in this
case the fact that the map $NS_m(X) \to NS_m(F)$ is surjective is
trivial. In fact, as we said, for a smooth hypersurface $F$ of odd
dimension, one has $dim\, NS_m(F)=1$. Previous Theorem \ref{main}
and Theorem \ref{mmain} give the following Corollary \ref{coro},
from which, taking into account (\ref{factor}), we obtain the
announced characterization for the factoriality of threefolds:

\medskip
\begin{corollary}
\label{coro} Let $X\subset \Ps^{2m+1}$ be a projective variety of
dimension $2m-1\geq 1$, with isolated singularities, complete
intersection of a smooth hypersurface $G$ of degree $k$, with a
hypersurface of degree $n$. Set $\lc =\mid \ic _{X,\Ps^{2m+1}}(n)
\mid $, and let $F\in\lc$ be a general hypersurface. If $n\geq k$
then $F$ is smooth. Moreover, if $n>max\{k, 2m+1\}$  then
$dim\,NS_m(X)=1$ if and only if $dim\,NS_m(F)=1$. In particular,
when $X$ is a threefold (i.e. $m=2$) then $X$ is factorial if and
only if $dim\, NS_2(F)=1$.
\end{corollary}
\medskip

Theorem \ref{mmain} should be compared with \cite{Lopez}, where
one describes the Picard group of a general complete intersection
surface containing a fixed smooth curve. Observe  also that,
assuming $k=1$, all previous results apply to any hypersurface
$X\subset \mathbb{P}^{2m}$ of degree $n$, with isolated
singularities.

The line of the proof of Theorem \ref{main} is the following. As
we said, $X$ only has isolated singularities by Proposition 4.3.6.
in \cite{Vogel}. Next, in order to prove that any projective
subvariety $S\subset X$ of dimension $m$ is homologous to a
multiple of the linear section $H^{m-1}_X$ of $X$, using
\cite{Fulton}, Example 15.3.2, we may assume that $S$ is integral
with isolated singularities. For such a subvariety, using our
hypotheses on $NS_m(F)$, the positivity assumption $n>k$ and a
suitable application of Hodge Index Theorem for $G$, we are able
to compare the double point formulae relative to the inclusions
$S\subset F$, $S\subset G$ and $S\cap R\subset X\cap R$, where $R$
is a general hypersurface of any degree $r\geq 1$ (so that $S\cap
R$ and $X\cap R$ are smooth). It turns out that $S\cap R$ is
homologous to a multiple $lH_{X\cap R}^{m-1}$ of the linear
section $H_{X\cap R}^{m-1}$ of $X\cap R$ in $H^{2m-2}(X\cap R;
\mathbb{C})$. To lift this homology to the whole $X$, we consider
a general hypersurface $R$ of degree $r\geq 1$, and a general
pencil $\rho: {\tilde{X}}\to \Ps^1$ of hypersurface sections
$\rho^{-1}(t)=X\cap R_t$ of $X$ of degree $r$ with $R=R_0$
(${\tilde{X}}=$ blowing-up of $X$ along the exceptional subset of
the pencil). The image $\tau(S-lH_{X}^{m-1})$ of the cycle
$S-lH_{X}^{m-1}$ through the Gysin morphism $\tau: H_{2m}({{X}};
\mathbb{C})\to H_{2m}({\tilde{X}}; \mathbb{C})$, maps to $0$ in
$H_{2m}(\tilde{X},\rho^{-1}(\cc); \mathbb{C})$ ($\cc=$ critical
locus of the pencil) because, using the Invariant Subspace Theorem
(see \cite{PS}, p.165-166), one may prove that
$H_{2m}(\tilde{X},\rho^{-1}(\cc); \mathbb{C})$ is canonically
embedded in $H^{2m-2}(X\cap R_t; \mathbb{C})$ (here we need that
the dimension of $X$ is odd). Therefore $\tau(S-lH_{X}^{m-1})$
comes from $H_{2m}(\rho^{-1}(\cc); \mathbb{C})$, and so it is $0$,
because by \cite{Dimca}, Theorem $(4.3)$, p. 161, we know that
this space is generated by the linear sections of the singular
fibres parametrized by $\cc$. It follows that also
$S-lH_{X}^{m-1}$ is $0$ (i.e. $S=lH_{X}^{m-1}$ in $H_{2m}(X;
\mathbb{C})$) because $\tau$ is injective.

More generally, the same argument we previously used to lift the
homology of algebraic cycles to $X$ applies to any cycle, i.e. one
has the following:
\medskip
\begin{proposition}
\label{lift} Let $X\subset \Ps^{N}$ be a complete intersection
projective variety of odd dimension $2m-1\geq 1$, with isolated
singularities.  Then for a general hypersurface $R$ of degree
$r\geq 1$ the Gysin morphism $H_{2m}(X; \mathbb{C})\to
H_{2m-2}(X\cap R; \mathbb{C})$ is injective.
\end{proposition}
\medskip

Observe that when $X$ is smooth, then previous Proposition
\ref{lift} follows from Lefschetz Hyperplane Theorem. We need
Proposition \ref{lift} for a comment on Theorem \ref{NL} below.

As for the proof of Theorem \ref{mmain}, first we prove that $F$
is smooth using a Bertini type of argument (which holds true also
when $X$ is of codimension two in $\Ps^{2m+2}$). Next we  show
that, in a certain sense, the classical Noether-Lefschetz argument
applies in our setting. More precisely, denote by $\qc\subset
\Ps^N$ the image of $\Ps^{2m+1}$ through the rational map
$\Ps^{2m+1} - -\to \Ps ^N$ defined choosing a basis of the linear
system $\lc$ ($N=dim \lc$). The variety $\qc$ only has isolated
singularities, and so one may regard $F$ as a general hyperplane
section of $\qc$, and may vary it in a general pencil $L\subset
{\Ps^N}^{*}$ of hyperplane sections. The monodromy action of the
pencil induces an orthogonal decomposition
\begin{equation}
\label{orth} H^{2m}(F; \mathbb{C})=I\oplus V,
\end{equation}
where $I$ is the  subspace of the invariant cocycles, and $V$ is
its orthogonal complement. Using a standard argument, we reduce
the proof of Theorem \ref{mmain} to prove that the monodromy of
the pencil irreducibly acts on $V$, i.e.  we prove the following:
\smallskip
\begin{theorem}
\label{Saito} The monodromy representation on $V$ for the family
of hypersurfaces of degree $n$ containing $X$ is irreducible.
\end{theorem}
\medskip
To prove this (see also Remark \ref{notasaito} below), using
\cite{La} and the theory of isolated singular points on complete
intersections as developed in \cite{Loj}, first we prove that $V$
is generated by the vanishing cocycles corresponding to the
hyperplane sections of $\qc$ which are {\it tangent} at some
regular point of $\qc$, and by a certain subspace of the space
generated by the vanishing cocycles defined by the
{\it{remaining}} singular hyperplane sections, i.e. by the
hyperplane sections of $\qc$ passing through its singularities
(except for the singularity of $\qc$ coming from the contraction
of $G$). Next we prove a basic lemma (i.e. Lemma \ref{trivial}
below) which states that the monodromy {\it trivially} acts on the
vanishing cocycles of the latter type. This implies that $V$ only
is  generated by the vanishing cocycles coming from tangential
sections, and so we may conclude the proof of Theorem \ref{Saito}
using the classical Zariski Theorem. To prove Lemma \ref{trivial},
using \cite{Loj}, first we reduce it to the case $X$ is a complete
intersection like $G\cap F'$, with $G$ a general hypersurface and
$F'$ a general hypersurface with a unique double point $q_1$ also
belonging to $G$. Then we conclude the proof of Lemma
\ref{trivial} by an \lq\lq ad hoc\rq\rq argument, relying on the
fact that one may realize the Milnor fibre of a general element
$F''$ of the linear system $\lc$ passing through $q_1$, as
contained in a {\it{fixed}} sphere which does not depend on $F''$
(see (\ref{varphi}) below). This argument does not apply to those
$F''$ corresponding to limit of tangential hyperplane sections of
$\qc$ for which there exists a sequence of regular contact points
converging to $q_1$. This case requires a separate analysis, which
in turn relies on the fact that, being $F'$ general, such $F''$
are parametrized by the dual variety of the tangent cone of $\qc$
at $q_1$, which is a nondegenerate and irreducible quadric in the
projective space parametrizing the hypersurfaces passing through
$q_1$.

Combining Theorem \ref{Saito} with \cite{Hu}, Corollary 1.1, we
obtain the following corollary:

\medskip
\begin{corollary}
\label{maincoro} Let $X\subset \Pcq$ be a projective threefold
with isolated singularities, complete intersection of a smooth
hypersurface  of degree $k$, with a hypersurface of degree $n$.
Assume $n>max\{k, 5\}$, set  $\lc =\mid \ic _{X,\Pcq}(n) \mid $,
and let $F\in\lc$ be a general hypersurface. Then Hodge conjecture
holds true for $F$.
\end{corollary}
\medskip

In fact, Theorem \ref{Saito} implies that $H^{2,2}(F)\subset I$.
Hence all the Hodge cycles of $F$ come from the Hodge cycles of a
desingularization of $\qc$, i.e. from some rational projective
complex manifold of dimension $5$, for  which, by \cite{Hu},
Corollary 1.1 (here we are forced to assume $X\subset \Pcq$),
Hodge conjecture holds true.

Both Theorem \ref{Saito} and Corollary \ref{maincoro} should be
compared with (\cite{OS}, Conjecture 0.2, Theorems 0.3 and 0.4,
and Corollary 0.5), where the authors prove  similar results, with
different assumptions. To this purpose, let we make two remarks.

First we notice one may prove that the subspace $I\subset
H^{2m}(F; \mathbb{C})$ defined by decomposition (\ref{orth}) is
the image of $H_{2m}(X; \mathbb{C})$ in $H_{2m}(F;
\mathbb{C})\simeq H^{2m}(F; \mathbb{C})$, and so, similarly as in
\cite{OS}, the subspace $V$ (for which our Theorem \ref{Saito}
states the irreducibility) is nothing but the quotient of
$H^{2m}(F; \mathbb{C})$ by the cycle classes coming from $X$. In
other words, with an analogous notation  as in \cite{OS}, one has
$V=H^{2m}(F; \mathbb{C})^{\text{van}}_{\perp X}$, and we may
restate previous Theorem \ref{Saito} as follows:
\smallskip
\begin{theorem}
\label{Saitodue}Let $X\subset \Ps^{2m+1}$ be a projective variety
of dimension $2m-1\geq 1$, with isolated singularities, complete
intersection of a smooth hypersurface of degree $k$, with a
hypersurface of degree $n>k$. Then the monodromy representation on
$H^{2m}(F; \mathbb{C})^{\text{van}}_{\perp X}$ for the family of
hypersurfaces $F$ of degree $n$ containing $X$ is irreducible.
\end{theorem}
\medskip

For the proof see {\it Proof of Theorem \ref{Saitodue}} in Section
3 below.

Next consider a projective surface $Z\subset \Ps^5$ whose ideal is
generated in degrees $\geq \delta$. Under mild assumptions on the
singularities of $Z$, one knows that $Z$ is contained in  smooth
hypersurfaces $G$ and $F$ of degree $k=\delta+1$ and $n=\delta+2$
(see \cite{OS}). From our Theorem \ref{main}  we also know that
$X=G\cap F$ only has isolated singularities, and therefore from
Corollary \ref{maincoro} (when $\delta >3$) the general
hypersurface $F$ of degree $n$ containing $X$ verifies Hodge
conjecture. A fortiori this holds true for a general hypersurface
of degree $n=\delta+2$ containing $Z$. So we see that, in the case
of a family of hypersurfaces of $\Ps^5$, our Corollary
\ref{maincoro} (at least when $\delta >3$, and for that concerns
the assertion on Hodge conjecture) implies Corollary 0.5 in
\cite{OS}.

As a further consequence of the proof of our Theorem \ref{Saito},
we may state a Noether-Lefschetz type Theorem for complete
intersections $\qc$ with isolated singularities, i.e. we are able
to prove the following:

\medskip
\begin{theorem}(Noether-Lefschetz Theorem with isolated
singularities) \label{NL} Let $\qc\subset \Ps^{N}$ be an
irreducible complete intersection projective variety with isolated
singularities, of odd dimension $2m+1\geq 3$. Assume that
$dim(NS_{m+1}(\qc))=1$. Then for any integer $r>>0$ and any
general hypersurface $R$ of degree $r$ one has $dim(NS_{m}(\qc\cap
R))=1$.
\end{theorem}
\medskip

For the proof see {\it Proof of Theorem \ref{NL}} in Section 3
below. Observe that, in view of Proposition \ref{lift}, the
assumption $dim(NS_{m+1}(\qc))=1$ in Theorem \ref{NL} is
necessary.

Finally we point out that, using  Theorem \ref{main}, we are able
to construct complete intersections $X$ of dimension $2m-1\geq 3$
with isolated singularities and with $dim(NS_m(X))=1$ (e.g.
factorial threefolds with isolated singularities), in some cases
also with many singularities: see Corollary \ref{examplesone} and
Corollary \ref{examplestwo} below. In particular we prove that the
asymptotic behavior of the maximal integer $r$ for which  there
exists a nodal factorial threefold in $\Ps^5$ complete
intersection of a smooth hypersurface of degree $n-1$ with a
hypersurface of degree $n$, with $|Sing(X)|=r$, is $n^5$. We also
stress that from \cite{Dimca}, Theorem (4.5), one may deduce that
for a nodal hypersurface $X\subset \Ps^{2m}$ of degree $n$ with at
most $m(n-2)$ nodes, one has $dim(NS_m(X))=1$.

Now we are going to prove the announced results.
\bigskip
\bigskip

\section{Proof of Theorem \ref{main} and consequences}

\begin{proof}[Proof of Theorem \ref{main}]
By Proposition 4.3.6. in \cite{Vogel} it follows that $X$ has at
most isolated singularities.

Now fix an integral subvariety $S\subset X$ of dimension $m$. To
prove Theorem \ref{main} it suffices to prove that $S$ is
homologous in $X$ to a multiple of the $m$-dimensional linear
section of $X$. To this aim first notice that, by \cite{Fulton},
Example 15.3.2, we may assume $Sing(S)\subset Sing(X)$, i.e. we
may assume  $S$ with isolated singularities. In particular, if
$R\subset \Ps^{2m+1}$ denotes a general hypersurface of degree
$r\geq 1$, then $C=S\cap R$ and $Y=X\cap R$ are smooth projective
varieties of dimensions $m-1$ and $2m-2$, with $C\subset Y$. In
$H^{2m-2}(Y; \mathbb{C})$ we may write
$$
C={\frac{d}{kn}}H_Y^{m-1}+\alpha,
$$
where $d$ is the degree of $S$, $H_Y$ is the general hyperplane
section of $Y$,  and $\alpha\in H^{2m-2}(Y; \mathbb{C})$ is a
primitive class, i.e. $\alpha.H_Y=0$ in $H^{2m}(Y; \mathbb{C})$.
We deduce:
$$
(C.C)_Y={\frac{d^2r}{kn}}+\alpha^2,
$$
where  $(C.C)_Y$ denotes  the self-intersection of $C$ in $Y$. On
the other hand, from the double point formula (\cite{Fulton}, p.
166), we know that
$$
(C.C)_Y=(c(i^*T_Y)c(T_C)^{-1})_{m-1},
$$
where $i$ denotes the inclusion $C\subset Y$, $T_Y$ and $T_C$ the
tangent bundles, and $c$  the total Chern class. Putting together
we obtain:
\begin{equation}
\label{alpa} \alpha^2=(c(i^*T_Y)c(T_C)^{-1})_{m-1}
-{\frac{d^2r}{kn}}.
\end{equation}

Besides  previous double point formula we may also consider the
two double point formulae corresponding to the inclusions
$S\subset F$ and $S\subset G$. More precisely, let $\pi: \Sigma
\to S$ be a desingularization of $S$. Denote by $f: \Sigma \to F$
and $g: \Sigma \to G$ the compositions of $\pi $ with the natural
inclusions. Following \cite{Fulton}, p. 166,  denote by
$\widetilde{\Sigma \times \Sigma}$ the blowing-up along the
diagonal, by $\df \subset \widetilde{\Sigma \times \Sigma}$ and
$\dg \subset \widetilde{\Sigma \times \Sigma}$ the {double point
schemes} of $f$ and $g$, by $D(f)\subset \Sigma $ and $D(g)\subset
\Sigma$ the {double point sets}, and by $\ddf \in A_0(D(f))$ and
$\ddg \in A_0(D(g))$  the {double point classes}. Applying the
double point formula  to $f$ and $g$ we obtain
\begin{equation}
\label{eq1} deg(\ddf) = (S.S)_F-(c(f^*T_F)c(T_{\Sigma})^{-1})_{m}
\end{equation}
and
\begin{equation}
\label{eq2} deg(\ddg) = (S.S)_G-(c(g^*T_G)c(T_{\Sigma})^{-1})_{m},
\end{equation}
where $(S.S)_F$ and $(S.S)_G$ represent  the self-intersection of
$S$ in $F$ and in $G$.

We claim that
\begin{equation}
\label{dd} deg(\ddf)=deg(\ddg).
\end{equation}
To prove this, denote  by $\varphi$ and $\gamma$ the natural maps
$\widetilde{\Sigma\times \Sigma}\to F\times F$ and
$\widetilde{\Sigma\times \Sigma}\to G\times G$, and by
$\Delta_F\subset F\times F$ and $\Delta_G\subset G\times G$ the
diagonals. Recall that $\df$ and $\dg$  are defined as the
residual schemes to the exceptional divisor  $E$ of
$\widetilde{\Sigma\times \Sigma}$, in $\varphi^{-1}(\Delta_F)$ and
in $\gamma^{-1}(\Delta_G)$. Since
$\varphi^{-1}(\Delta_F)=\gamma^{-1}(\Delta_G)$, then we have
\begin{equation}
\label{ugu} \df = \dg.
\end{equation}
From (\ref{ugu}) and  (\cite{Fulton}, p. 166), it follows that to
prove (\ref{dd}) it suffices to show that the residual
intersection classes $\dtf$ and $\dtg $ coincide in $A_0(\df
)=A_0(\dg )$. To this purpose, notice that by (\cite{Fulton},
Theorem 9.2) we have
\begin{equation}
\label{novef} \dtf =\{c(N_{\Delta_F}\otimes \oc (-E))\cap s(\df ,
\widetilde{\Sigma\times \Sigma})\}_0
\end{equation}
and
\begin{equation}
\label{noveg} \dtg =\{c(N_{\Delta_G}\otimes \oc (-E))\cap s(\dg ,
\widetilde{\Sigma\times \Sigma})\}_0,
\end{equation}
where $c$ and $s$ denote Chern and Segre classes, and $\oc (-E)$,
$N_{\Delta_F}$ and $N_{\Delta_G}$ are the pull-back on
$\varphi^{-1}(\Delta_F)=\gamma^{-1}(\Delta_G)$ of
$\mathcal{O}_{\widetilde{\Sigma\times \Sigma}}(-E)$ and of the
normal bundles of ${\Delta_F}$ and ${\Delta_G}$ in $F\times F$ and
$G\times G$. Since these normal bundles are isomorphic to the
tangent bundles of $F$ and $G$ and the Chern polynomials of both
$F$ and $G$ only depend  on the hyperplane class, and since $\pi $
is an isomorphism outside of a finite set of $S$, then both
$c(N_{\Delta_F})$ and $c(N_{\Delta_G})$ are the identity in $A(\df
)=A(\dg )$. In particular $c(N_{\Delta_F})=c(N_{\Delta_G})$ in
$A(\df )=A(\dg )$. Therefore from (\ref{ugu}), (\ref{novef}) and
(\ref{noveg}), we obtain (\ref{dd}).

Now we notice that our assumption on $NS_m(F)$ implies that
$(S.S)_F =\frac{d^2}{n}$. On the other hand, by Hodge Index
Theorem for $G$ (see \cite{Hartshorne}, Theorem 5.2, pg. 435), we
have $(-1)^m(S-\frac{d}{k}H_G^m)^2\geq 0$ on $G$ ($H_G=$ general
hyperplane section of $G$), and so $(-1)^m((S.S)_G
-\frac{d^2}{k})\geq 0$. Comparing with (\ref{alpa}), (\ref{eq1}),
(\ref{eq2}) and (\ref{dd}),  and taking into account that $k<n$,
we obtain that $(-1)^{m-1}\alpha^2$ is less than or equal to
\begin{equation}
\label{eq3}(-1)^{m-1}\left[(c(i^*T_Y)c(T_C)^{-1})_{m-1}-\frac{r}{k-n}
\left[(c(f^*T_F)-c(g^*T_G))c(T_{\Sigma})^{-1}\right]_{m}\right].
\end{equation}
Using the exact sequences $0\to T_C\to T_{\Sigma}|_{C}\to
\oc_C(r)\to 0$, $0\to T_Y\to T_{F}|_{Y}\to \oc_Y(r)\oplus
\oc_Y(k)\to 0$ and $0\to T_Y\to T_{G}|_{Y}\to \oc_Y(r)\oplus
\oc_Y(n)\to 0$, one may compare the intersection numbers appearing
in (\ref{eq3}). It turns out that the number in the formula
(\ref{eq3}) is $0$, therefore we have:
$$
(-1)^{m-1}\alpha^2\leq 0.
$$
By Hodge Index Theorem for $Y$ it follows  that $\alpha=0$. In
other words, for any integer $r\geq 1$ and for a general
hypersurface $R\subset\Ps^{2m+1}$ of degree $r$, one has
$$
S\cap R= \frac{d}{kn}H_{X\cap R}^{m-1}\quad{\text{in}}\quad
H^{2m-2}(X\cap R; \mathbb{C})
$$
($H_{X\cap R}=$ general hyperplane section of $X\cap R$). At this
point, to conclude the proof of Theorem \ref{main} it suffices to
prove Proposition \ref{lift}.

\bigskip \noindent {\it Proof of Proposition \ref{lift}}. Consider
a general pencil $\{X\cap R_t\}_{t\in\Ps^1}$ of
hypersurface sections of $X$, with $deg(R_t)=r$. We may regard the
pencil as the set of fibres of a projective morphism
$$
\rho: {\tilde{X}}\to \Ps^1,
$$
where ${\tilde{X}}$ is the blowing-up of $X$ along the exceptional
subset. Let $\cc\subset\Ps^1$ be the critical locus of $\rho$, put
$U=\Ps^1-\cc$, and consider the natural exact sequence
\begin{equation}
\label{homoseq} H_{2m}(\rho^{-1}(\cc); \mathbb{C})\to
H_{2m}(\tilde{X}; \mathbb{C})\to H_{2m}(\tilde{X},\rho^{-1}(\cc);
\mathbb{C}).
\end{equation}
Applying Lefschetz Duality to the pair
$(\tilde{X},\rho^{-1}(\cc))$ (\cite{Spanier}, p. 297), we obtain a
natural isomorphism
\begin{equation}
\label{iso1}H_{2m}(\tilde{X},\rho^{-1}(\cc); \mathbb{C})\simeq
H^{2m-2}(\rho^{-1}(U); \mathbb{C}).
\end{equation}
The Leray spectral sequence of the  restriction $\rho^{-1}(U)\to
U$ collapses in the term $E_2$ (see \cite{PS}, p. 166), and
therefore we have an isomorphism:
\begin{equation}
\label{iso2} H^{2m-2}(\rho^{-1}(U); \mathbb{C})\simeq
\oplus_{j\geq 0}\,H^j(U,R^{2m-2-j}\rho_{*} \mathbb{C}).
\end{equation}
Since for any $t\in U$ the fiber $\rho^{-1}(t)=X\cap R_t$ is a
smooth projective complete intersection of {\it{even}} complex
dimension $2m-2$, it follows that $H^j(U,R^{2m-2-j}\rho_{*}
\mathbb{C})=0$ for any $j\geq 1$ odd. For the same reason, when
$j\geq 1$ is even, the local system $R^{2m-2-j}\rho_{*}
\mathbb{C}$ on $U$ has rank $1$ and has a global section with no
zeros, corresponding to the linear section of $X$ with a general
subspace of $\Ps^{2m+1}$ of codimension $(2m-2-j)/2$. Hence
$R^{2m-2-j}\rho_{*} \mathbb{C}$ is isomorphic to the constant
sheaf $\mathbb{C}$. So we have $H^j(U,R^{2m-2-j}\rho_{*}
\mathbb{C})=H^j(U; \mathbb{C})$, which again vanishes when $j\geq
1$ is even, by Lefschetz Duality (we may assume that $\cc$ is non
empty). Therefore, from (\ref{iso2}), we deduce that the natural
map
$$
H^{2m-2}(\rho^{-1}(U); \mathbb{C})\to H^0(U,R^{2m-2}\rho_{*}
\mathbb{C})
$$
is an isomorphism. Taking into account that, for any $t\in U$,
$H^0(U,R^{2m-2}\rho_{*} \mathbb{C})$ identifies with the invariant
subspace $H^{2m-2}(\rho^{-1}(t);\mathbb{C})^{\text{inv}}$, it
follows a natural inclusion:
\begin{equation}
\label{inclusion} H^{2m-2}(\rho^{-1}(U); \mathbb{C})\subset
H^{2m-2}(\rho^{-1}(t);\mathbb{C}).
\end{equation}

Now fix a cycle $a\in H_{2m}({X}; \mathbb{C})$ which restricts to
$0$ in $H_{2m-2}({X\cap R_t}; \mathbb{C})\simeq
H^{2m-2}(\rho^{-1}(t); \mathbb{C})$ ($t\in U$), and let $\tau(a)$
be the image of $a$ through the Gysin morphism $\tau:
H_{2m}({{X}}; \mathbb{C})\to H_{2m}({\tilde{X}}; \mathbb{C})$ (see
\cite{Fulton}, Example 19.2.1). Using (\ref{homoseq}),
(\ref{iso1}) and (\ref{inclusion}) we see that the map
$H_{2m}(\tilde{X}; \mathbb{C})\to H_{2m}(\tilde{X},\rho^{-1}(\cc);
\mathbb{C})$ sends $\tau(a)$ to $0$. From the exact sequence
(\ref{homoseq}) it follows that $\tau(a)$ comes from
$H_{2m}(\rho^{-1}(\cc); \mathbb{C})$. This space is the direct sum
of a finite number of spaces like $H_{2m}(X\cap R_{t'};
\mathbb{C})$, where $X\cap R_{t'}$ is a projective complete
intersection of dimension $2m-2$ with isolated singularities. From
\cite{Dimca}, Theorem $(4.3)$, p. 161, we know that such a space
has dimension $1$. Hence $\tau(a)$ is equal to a certain multiple
of $\tau({H}_{X}^{m-1})$ in $H_{2m}(\tilde{X}; \mathbb{C})$. Since
$a$ restricts to $0$ in $H_{2m-2}({X\cap R_t}; \mathbb{C})$ it
follows that $\tau(a)=0$, which  implies that $a=0$ because $\tau$
is injective.

This concludes the proof of Proposition \ref{lift} and, as we said
before, the proof  of Theorem \ref{main}.
\end{proof}

\bigskip

One may consider Theorem \ref{main} as a method to construct
complete intersections $X$ of dimension $2m-1\geq 3$, with
isolated singularities and with $dim(NS_m(X))=1$ (e.g. factorial
threefolds with isolated singularities). A  question arising
naturally from this remark is how many (and what kind of)
singularities one may produce in this fashion. To this purpose we
are able to prove the following Corollary \ref{examplesone} and
Corollary \ref{examplestwo}. First we need  the following:

\medskip
\begin{lemma}
\label{nodes} Let $F\subset \Ps^b$ ($b\geq 2$) be a projective
smooth hypersurface. Fix an integer $k>0$, and  $r$ points $\Sigma
=\{p_1,\dots,p_r\}$ on $F$. Assume that
$r<\left(\frac{b+1+k}{b+3}\right)^b$, and that the ideal of
$\Sigma$ is generated in degree $\leq \frac{b+1+k}{b+3}$. Then
there exists a projective smooth hypersurface $G\subset \Ps^b$ of
degree $k$ such that the singular locus of the complete
intersection $X=G\cap F$ is $\Sigma$, and each point $p_i$ is an
ordinary double point for $X$.
\end{lemma}

\begin{proof}
Denote by $\pi:\Ps_1\to\Ps^b$ the blowing-up  of $\Ps^b$ along
$\Sigma$,  by $E=\sum_{i=1}^{r} E_i$ its exceptional divisor, and
by $\tilde H_i$ the strict transform  in $\Ps_1$ of the tangent
hyperplane of $F$ at $p_i$. Put $H_i=\tilde H_i\cap E_i$, denote
by $\rho:\Ps_2\to\Ps_1$ the blowing-up  of $\Ps_1$ along
$\sum_{i=1}^{r} H_i$, and by $L=\sum_{i=1}^{r} L_i$ its
exceptional divisor. Consider the linear system:
$$
\mid D\mid = \mid \rho^*(k\pi^*(H)-E)-L\mid,
$$
where $H\subset \Ps^b$ denotes a hyperplane divisor. As a first
step we prove that $\mid D\mid$ is base-point free.

To this purpose first we prove that $\mid D\mid$ is base-point
free on $L$. To this aim notice that $H_i\simeq\Ps^{b-2}$, and
that its normal bundle in $\Ps_1$ is isomorphic to
$\oc_{H_i}(-1)\oplus \oc_{H_i}(1)$. Therefore the $\Ps^1$-bundle
$L_i\to H_i$ is isomorphic to $\Ps(\oc_{H_i}\oplus \oc_{H_i}(2))$.
The divisor $D$ restricts on each component  $L_i$ to a moving
section  of $L_i\to H_i$, which defines a base-point free linear
system on $L_i$. Hence, by the defining sequence of $\oc _{L}(D)$,
it follows that to prove $\mid D\mid$ is base-point free on $L$ it
suffices to prove $h^1(\Ps_2, \oc_{\Ps_2}(D-L))=0$. On the other
hand, since
$D-L=K_{\Ps_2}+(b+1+k)(\pi\circ\rho)^*(H)-b\rho^*(E)-3L$
($K_{\Ps_2}$= canonical divisor of $\Ps_2$), then by
Kawamata-Vieweg Theorem it suffices to prove that the divisor
$A=(b+1+k)(\pi\circ\rho)^*(H)-b\rho^*(E)-3L$ is big and nef.

Since the restrictions of $A$  to the strict transform $\tilde E$
of $E$ and to $L$ are very ample, to prove that $A$ is nef it
suffices to prove that $A.C_2\geq 0$ for the strict transform
$C_2\subset \Ps_2$ of any irreducible curve  $C\subset\Ps^b$.
Taking into account that $\rho^*(E).C_2\geq L.C_2$, and using the
projection formula, we see that $A.C_2\geq
((b+1+k)\pi^*(H)-(b+3)E).C_1$ ($C_1$= strict transform of $C$ in
$\Ps_1$), which is $\geq 0$ for our assumption on the degree of
the generators of the ideal of $\Sigma$. This proves that $A$ is
nef.

To prove that $A$ is big, first we notice that, using previous
description of the bundle $L_i\to H_i$ and the intersection
formulae for a blowing-up appearing in \cite{Fulton}, p. 67, one
may prove that ${L_i}^j.\rho^*(E_i)^{b-j}=1$ for $0\leq j\leq b$
even, and $0$ otherwise. It follows that $A^b\geq
\left[(b+1+k)^b-(b+3)^b\right]r$, which is $>0$ for our assumption
on $r$. This proves that $A$ is big, and so that $\mid D\mid$ is
base-point free on $L$.

A similar computation proves that $\mid D\mid$ is base-point free
on $\tilde E$, and out of $\tilde E+L$. This proves that $\mid
D\mid$ is base-point free. Therefore, by Bertini Theorem, any
general divisor $D\in \mid D\mid$ is smooth. Then one may choose
$G=(\pi\circ\rho)_*(D)$, because similar argument and computation
as before prove that $G$ is isomorphic to $D$,  that $Sing(G\cap
F)=\Sigma$, and that each $p_i$ is a node for $G\cap F$.
\end{proof}
\medskip

We are in position to prove the announced corollaries:

\medskip
\begin{corollary}
\label{examplesone} Fix integers $m\geq 2$, $k$ and $r$, put
$q=\frac {2m+2+k} {2m+4}$, and assume that $r\leq {\binom
{q}{2m+1}}$. Then for any $n >k$ there exists an integral
projective complete intersection $X\subset \Ps^{2m+1}$ of
bi-degree $(k,n)$ whose singular locus consists of exactly $r$
ordinary double points, and such that $dim(NS_m(X))=1$.
\end{corollary}

\begin{proof} Fix $r$  general points $p_1,\dots,p_r$ in
$\Ps^{2m+1}$. Since $r< {\binom {n+2m+1}{2m+1}}$ then there exists
a general hypersurface $F$ of degree $n$ passing through
$p_1,\dots,p_r$. By Noether-Lefschetz Theorem we have
$dim(NS_m(F))=1$. On the other hand, combining our assumption
$r\leq {\binom {q}{2m+1}}$ with Lemma \ref{nodes}  and \cite{GM},
Corollary 1.6, we deduce the existence of a smooth hypersurface
$G$ of degree $k$ such that the singular locus of the complete
intersection $X=G\cap F$ consists of exactly $r$ nodes at the
points $p_1,\dots,p_r$. By our Theorem \ref{main}, one has
$dim(NS_m(X))=1$.
\end{proof}

Next result states  that  the examples exhibited by previous
Corollary \ref{examplesone} in $\Ps^5$ are asymptotically sharp,
i.e. one has:

\medskip
\begin{corollary}
\label{examplestwo} For any integer $n\geq 3$ denote by $\nu(n)$
the maximal integer $r$ for which there exists a nodal factorial
threefold in $\Ps^5$ complete intersection of a smooth
hypersurface of degree $n-1$ with a hypersurface of degree $n$,
with $|Sing(X)|=r$. Then there exist positive constants
$\gamma_1>0$ and $\gamma_2>0$ such that
$$
\gamma_1 \leq \frac {\nu(n)} {n^5}\leq \gamma_2
$$
for any $n\geq 3$.
\end{corollary}

\begin{proof} From Corollary \ref{examplesone} we deduce that $\frac {\nu(n)}
{n^5}$ is bounded from below by some positive constant for $n>>0$.
On the other hand, by \cite{Cheltsov2} and \cite{CD2} we know that
$\nu(n)\geq 1$ for any $n\geq 3$. This proves the existence of
$\gamma_1$. As for $\gamma_2$, recall that the defect of a nodal
factorial threefold $X$ vanishes. If such a threefold $X$ is a
complete intersection in $\Ps^5$ of a smooth hypersurface of
degree $n-1$ with a hypersurface of degree $n$, this means that
the value of the Hilbert function of the singular locus of $X$ at
level $3n-7$ is $|Sing(X)|$ (see \cite{Cynk}). Therefore one has
$|Sing(X)|\leq \binom {3n-2}{5}$.
\end{proof}

\bigskip

Our numerical assumption in the proof of Lemma \ref{nodes}
certainly is not the best possible. It only  is of the simplest
form we were able to conceive. We also notice that the proof of
Lemma \ref{nodes} can be generalized to worse singularities, and
that one may state a similar result as in Corollary
\ref{examplestwo} for threefolds complete intersections in $\Ps^5$
of bi-degree $(k,n)$ ($k<n$) with $k$ not too far from $n-1$. We
decided not to push here this investigation further. We have in
mind to give more information on this subject in a forthcoming
paper.

\bigskip
\bigskip

\section{Proof of Theorem \ref{mmain} and consequences}

\begin{proof}[Proof of Theorem \ref{mmain}] For any $p\in Sing(X)$, let
$\lc_p\subset \lc$ be the closed set of all $F\in \lc$ such that
$p\in Sing(F)$. By Bertini Theorem, the singular locus of any
general $F\in \lc$ is contained in $X$. Since $X$ is the complete
intersection of $F$ and $G$, then any singular point of $F$ has to
be also a singular point for $X$. In other words, for a general
$F\in\lc$ one has $Sing(F)\subset Sing(X)$. Therefore, in order to
prove that the general $F\in \lc$ is smooth, it suffices to prove
that, for any $p\in Sing(X)$, $\lc- \lc_p$ is  non empty. To this
purpose, fix $D\in \lc_p$, and denote by $R$ a hypersurface of
degree $n-k$ such that $p\notin R$. Let $D_1\in \lc$ be the
hypersurface defined by the equation $d+rg=0$, where $g=0$, $d=0$
and $r=0$ are the equations defining $G$, $D$ and $R$ (when $n=k$,
define $D_1$ by the equation $d+g=0$). Taking into account that
$G$ is smooth, computing derivatives one sees that $D_1$ is smooth
at $p$, and so $\lc- \lc_p$ is a non empty subset of $\lc$.

Next, we are going to prove that if $n>max\{k,2m+1\}$ then, for a
general $F\in \lc$, the map $NS_m(X)\to NS_m(F)$ is surjective.
The proof is an adaptation of the classical Noether-Lefschetz
argument.

Choosing a basis of the linear system $\lc$, we may define a
rational map:
\begin{equation}
\label{rationalmap} \Ps^{2m+1} --\to \Ps:=\Ps ^N
\end{equation}
($N=dim \lc$), whose resolution is represented by the blowing-up
$\pc$ of $\Ps^{2m+1}$ along $X$, equipped with natural maps $\pc
\to \Ps^{2m+1}$ and $\pc \to \Ps$ (see \cite{Hartshorne}, p. 168).
By \cite{Fulton}, p. 437, B.6.10., we know that locally $\pc$ is
the hypersurface of $\mathbb{A}^{2m+1} \times \Ps^1$ defined by
the equation
\begin{equation}
\label{blow} gu_0-fu_1=0,
\end{equation}
where $g=0$ and $f=0$ are the local equations of $G$ and $F$ in
$\mathbb{A}^{2m+1}$, and $u_0$, $u_1$ are coordinates in $\Ps^1$.
Denote by $\qc\subset \Ps$ the image of $\pc$. The map $\pc \to
\Ps$ sends all the points of the strict transform $\tG$ of $G$
(and only them) to a singular point $q_{\infty}$ of $\qc$, and,
since $n>k$, the global sections of the pull-back of
$\oc_{\Ps}(1)$ separate points and tangent vectors of $\pc$ out of
$\tG$, i.e. the map $\pc \to \Ps$ induces an algebraic isomorphism
\begin{equation}
\label{gtilde} \pc - \tG \simeq \qc - \{q_{\infty}\}.
\end{equation}
From (\ref{blow}) and (\ref{gtilde}) one deduces that, besides the
point $q_{\infty}$, the singular locus of  $\qc$ consists of a
certain finite number of points $q_1,\dots,q_r$ corresponding to
the singular points of $X$. Hence we have
$Sing(\qc)=\{q_1,\dots,q_r, q_{\infty}\}$, and
$Sing(\pc)=\{q_1,\dots,q_r\}$ (actually each point $q_i\neq
q_{\infty}$ is a double point, and it is a node if and only if the
corresponding singular point of $X$ is).

For any $x$ in the dual space $\Ps ^*$, denote by $H_x\subset \Ps$
the corresponding hyperplane, and by $F_x$ the corresponding
hypersurface of $\lc $. Next denote by $\dc \subset \Ps ^*$ the
discriminant variety of $\lc $, i.e. the variety parametrizing the
singular hypersurfaces of $\lc $. Such a variety has $r+2$
components: the dual variety $\qc ^*$ of $\qc$, the hyperplane
$\hc_{\infty}$ corresponding to the singular point $q_{\infty}$
(i.e. corresponding to the reducible hypersurfaces of $\lc$
containing $G$), and the $r$ hyperplanes $\hc_1,\dots,\hc_r$,
corresponding to the remaining singular points $q_1,\dots,q_r$.

Now fix a general point $t\in \Ps ^*$ and let $L$ be a general
line through $t$. This line meets each hyperplane $\hc_i$, $i\in
\{1,\dots,r,{\infty}\}$, in a certain point $a_i$, and meets $\qc
^*$ transversally in certain smooth points  $a_{r+1}, \dots ,a_s$.
All the points $a_{i}$, $i\neq \infty$,  correspond to irreducible
hypersurfaces $F_{a_{i}}$   in $\lc$ with a unique double point.
When $r+1\leq i\leq s$, then $F_{a_{i}}$ corresponds to a tangent
hyperplane section of $\qc$, and therefore its unique double point
is ordinary. The point $a_{\infty}$ corresponds to a reducible
hypersurface $M\cup G$ containing $G$:
\begin{equation}
\label{mpiug} F_{a_{\infty}}=M\cup G.
\end{equation}
The intersection
\begin{equation}
\label{baselocus} B=H_t\cap H_{a_{\infty}}\cap \qc
\end{equation}
defines the exceptional subset of  $\qc$ with respect to $L$. By
(\ref{gtilde}), we may regard it as a subset of $\pc$. Denote by
$\tp $ the blowing-up of $\pc$ along $B$, and by $\rc$ a
desingularization of $\tp $. Set $f: \rc \to L$ the natural
projection. The restriction
$$
f: \rc -f^{-1}(\{ a_1, \dots ,a_s ,a_{\infty}\}) \to  L-\{ a_1,
\dots ,a_s ,a_{\infty}\}=L-\dc
$$
is a smooth proper map. Hence the fundamental group $\pi_1(L-\dc,
t)$ acts by monodromy on $f^{-1}(t)\simeq F_t$, and so on
$H^{2m}(F_t; \mathbb{C})$. By the Invariant Subspace Theorem
\cite{PS}, p. 165-167, we know that there is an orthogonal
decomposition:
$$
H^{2m}(F_t; \mathbb{C})=I\oplus V,
$$
where $I$ is the  subspace of the invariant cocycles, and $V$ is
its orthogonal complement. If $j$ denotes the natural inclusion
$F_t\subset \rc$, then we also have $I=j^*H^{2m} (\rc ;
\mathbb{C})$ from which, using Poincar\`e duality, we get
$$
V=Ker(H^{2m}(F_t ; \mathbb{C})\to H^{2m+2} (\rc ; \mathbb{C}))
\simeq Ker(H_{2m}(F_t ; \mathbb{C})\to H_{2m}(\rc ; \mathbb{C})).
$$

We notice that
\begin{equation}
\label{VeKer} Ker(H_{2m}(F_t ; \mathbb{C})\to H_{2m}(\rc
;\mathbb{C}))=Ker(H_{2m}(F_t ; \mathbb{C})\to
H_{2m}(\rc-f^{-1}(a_{\infty}) ;\mathbb{C})).
\end{equation}
In fact, to prove (\ref{VeKer}), taking into account that
$F_t\subset \rc-f^{-1}(a_{\infty})$, it suffices to prove that the
natural map $H_{2m}(\rc-f^{-1}(a_{\infty}) ; \mathbb{C})\to
H_{2m}(\rc ;\mathbb{C})$ is injective. By Lefschetz and Poincar\`e
dualities we have natural isomorphisms
$H_{2m}(\rc-f^{-1}(a_{\infty}) ;\mathbb{C})\simeq
H^{2m+2}(\rc,f^{-1}(a_{\infty});\mathbb{C})$ and $H_{2m}
(\rc;\mathbb{C})\simeq H^{2m+2}(\rc;\mathbb{C})$. So, to prove
(\ref{VeKer}), it suffices to prove that the natural map $H^{2m+2}
(\rc,f^{-1}(a_{\infty});\mathbb{C})\to H^{2m+2}(\rc;\mathbb{C})$
is injective. This follows from the vanishing of
$H^{2m+1}(f^{-1}(a_{\infty}) ;\mathbb{C})$. To prove this, using
(\ref{blow}), first one sees that the strict transform $\tG$ of
$G$ in $\pc$  does not meet the singular locus of $\pc$. It
follows that $f^{-1}(a_{\infty})$ simply is the union $G_1\cup
M_1$ of the strict transforms of $G$ and $M$  in $\rc$ (compare
with (\ref{mpiug})). Since  $G_1\simeq G$ and $M_1$ is isomorphic
to the blowing-up of the general hypersurface $M\subset
\Ps^{2m+1}$ of degree $n-k$ along the smooth complete intersection
$M\cap X$ of dimension $2m-2$, and $G_1\cap M_1\simeq G\cap M$,
then one may compute the terms of the Mayer-Vietoris sequence of
the pair $(G_1,M_1)$ which allow to control
$H^{2m+1}(f^{-1}(a_{\infty}) ;\mathbb{C})$. It turns out that
$H^{2m+1}(f^{-1}(a_{\infty}) ;\mathbb{C})=0$.

From (\ref{VeKer}) and the homology exact sequence of the pair
$(\rc - f^{-1}(a_{\infty}), F_t)$ we deduce
\begin{equation}
\label{deco} V\simeq Im(H_{2m+1}(\rc-f^{-1}(a_{\infty}),F_t ;
\mathbb{C})\to H_{2m}(F_t ;\mathbb{C})).
\end{equation}
Now, as in \cite{La}, pg. 35, Fig. 1, for any $1\leq i \leq s$ fix
a closed disk $\Delta_i\subset L-\{a_{\infty}\} \simeq \bC$ with
center $a_i$ and radius $0<\rho<<1$, and, in
$\bC-\bigcup_{i=1}^{s}\Delta_i^{\circ}$ ($\Delta_i^{\circ}=$
interior of $\Delta_i$), choose a $C^{\infty}$ path $l_i$ from $t$
to $a_i+\rho$ with no self-intersection points and such that
$l_i\cap l_j=\{t\}$ for $i\neq j$. Using the same argument as in
\cite{La}, (5.3.1) and (5.3.2), one may prove a direct
decomposition
$$
H_{2m+1}(\rc-f^{-1}(a_{\infty}),F_t ;
\mathbb{C})\simeq{\oplus}_{i=1}^s H_{2m+1}
(f^{-1}(\Delta_i),f^{-1}({a_i+\rho}) ; \mathbb{C}).
$$
If we denote by $V_i$ the image of each $H_{2m+1}
(f^{-1}(\Delta_i),f^{-1}({a_i+\rho});\mathbb{C})$ in $H_{2m}(F_t
;\mathbb{C})$ $\simeq$  $H^{2m}(F_t ;\mathbb{C})$, then by
(\ref{deco}) we get a decomposition:
\begin{equation}
\label{dec} V=V_1+\dots+V_s.
\end{equation}
Notice that each path $l_i$ induces a $C^{\infty}$- diffeomorphism
$f^{-1}({a_i+\rho})\simeq F_t$, and so an isomorphism only
depending on $l_i$:
\begin{equation}
\label{diff} H_{2m}(f^{-1}({a_i+\rho}) ; \mathbb{C})\simeq H_{2m}
(F_t ;\mathbb{C}),
\end{equation}
which in turn identifies
\begin{equation}
\label{ident} Im(H_{2m+1}(f^{-1}(\Delta_i),f^{-1}({a_i+\rho});
\mathbb{C})\to H_{2m}(f^{-1}({a_i+\rho}) ; \mathbb{C}))\simeq V_i.
\end{equation}

When $r+1\leq i\leq s$, we recognize in $V_i\subset H^{2m}(F_t
;\mathbb{C})$ the subspace generated by the \lq\lq classical\rq\rq
vanishing cocycle corresponding to a tangent hyperplane section of
$\qc$ (see \cite{La}, \cite{Voisin}). For the remaining subspaces
we claim that
\begin{equation}
\label{van} V_i=0 \quad {\text{for any}}\, 1\leq i\leq r.
\end{equation}
To prove (\ref{van}), fix an index $1\leq i\leq r$ and denote by
$g$ the natural projection $\tp\to L$, so that $f$ is the
composition of $g$ with the desingularization $\rc\to\tp$. By
\cite{Loj}, p. 28, we know that near to the isolated singular
point $q_i\in\tp$, the pencil $g: \tp\to L$ defines a Milnor
fibration with Milnor fiber
\begin{equation}
\label{Milnorfiber} g^{-1}(a_i+\rho)\cap D_i,
\end{equation}
where $D_i$ denotes a closed ball of the ambient space in which
$\tp$ is embedded, with center $q_i$ and positive small radius
$\epsilon$ with $\rho << \epsilon<< 1$. Set
\begin{equation}
\label{defofI} I_i=Im(H_{2m}(g^{-1}(a_i+\rho)\cap D_i;
\mathbb{C})\to H_{2m}(g^{-1}(a_i+\rho);\mathbb{C})).
\end{equation}
Observe that $g^{-1}(a_i+\rho)$ is canonically isomorphic to
$f^{-1}(a_i+\rho)$. Hence, via $l_i$,  by (\ref{diff}) we may
regard $V_i$ and $I_i$ both contained in
$H_{2m}(g^{-1}(a_i+\rho);\mathbb{C})$ and in $H_{2m}(F_t
;\mathbb{C})\simeq H^{2m}(F_t ;\mathbb{C})$.

Since $g^{-1}(\Delta_i)- D^{\circ}_i\to \Delta_i$ is a trivial
fibre bundle ($D^{\circ}_i$= interior of $D_i$), using Excision
Axiom and Leray-Hirsch Theorem (\cite{Spanier}, p. 200 and 258),
one sees that the inclusion $(g^{-1}(a), g^{-1}(a)\cap D_i)\subset
(g^{-1}(\Delta_i),g^{-1}(\Delta_i)\cap D_i)$ induces natural
isomorphisms:
\begin{equation}
\label{foranya} H_{2m}(g^{-1}(a), g^{-1}(a)\cap D_i
;\mathbb{C})\simeq H_{2m}(g^{-1}(\Delta_i),g^{-1}(\Delta_i)\cap
D_i ;\mathbb{C})
\end{equation}
for any $a\in \Delta_i$. From (\ref{defofI}), (\ref{foranya}), the
homology sequence of the pair $(g^{-1}(a_i+\rho),
g^{-1}(a_i+\rho)\cap D_i)$,  and the conic structure of
$g^{-1}(\Delta_i)\cap D_i$ (\cite{Loj}, Lemma (2.10)), which
implies that $H_{2m}(g^{-1}(\Delta_i),g^{-1}(\Delta_i)\cap D_i
;\mathbb{C})\simeq H_{2m}(g^{-1}(\Delta_i);\mathbb{C})$, it
follows the natural exact sequence:
\begin{equation}
\label{nexact} 0\to I_i\to H_{2m}(g^{-1}(a_i+\rho);\mathbb{C})\to
H_{2m}(g^{-1}(\Delta_i);\mathbb{C}).
\end{equation}
Since we may regard the inclusion $g^{-1}(a_i+\rho)\subset
g^{-1}(\Delta_i)$ as the composition of the isomorphism
$g^{-1}(a_i+\rho)\simeq f^{-1}(a_i+\rho)$ with the inclusion
$f^{-1}(a_i+\rho)\subset f^{-1}(\Delta_i)$, followed by the
desingularization $f^{-1}(\Delta_i)\to g^{-1}(\Delta_i)$, from
(\ref{ident}) it follows that
$$
V_i\subset Ker(H_{2m}(g^{-1}(a_i+\rho);\mathbb{C})\to
H_{2m}(g^{-1}(\Delta_i);\mathbb{C}))
$$
and therefore, from (\ref{nexact}), we obtain
\begin{equation}
\label{VI} V_i\subset I_i.
\end{equation}

\medskip
\begin{remark} \label{milnornumber}
Notice that from the local description (\ref{blow}), it follows
that  the singularities of $\pc$, and hence of $\tp$, are all
locally complete intersection isolated singularities. Then the
Milnor fiber $g^{-1}(a_i+\rho)\cap D_i$ defined by the pencil $g$
around $q_i$ is the Milnor fibre of the isolated complete
intersection singularity $(g^{-1}(a_i)\cap D_i, q_i)$. Therefore
$g^{-1}(a_i+\rho)\cap D_i$ has the homotopy type of a bouquet of
$2m-$spheres  contained in $g^{-1}(a_i+\rho)$, and these
$2m-$spheres, as cycle classes, span $I_i$. In particular
$H_{2m-1}(g^{-1}(a_i+\rho)\cap D_i ;\mathbb{C})=0$, which implies
that the right map in (\ref{nexact}) actually is surjective
(compare with \cite{Loj}, pp. 7, 73-76, 121). The number $\mu$ of
$2m-$spheres occurring in the bouquet is called the Milnor number
of the singularity $(g^{-1}(a_i)\cap D_i, q_i)$.
\end{remark}
\medskip

At this point we need the following basic lemma.

\smallskip
\begin{lemma}
\label{trivial} For any $1\leq i\leq r$, the  group $\pi_1(L-\dc,
t)$ trivially  acts on $I_i$.
\end{lemma}
\smallskip

\begin{proof}[Proof of Lemma \ref{trivial}]
Our first step consists in proving that one may assume  $X$ with a
unique ordinary double point.

To this purpose, fix an integer $i\in \{1,\dots,r\}$, i.e. fix a
singular point $q_i$ of $\qc -\{q_{\infty}\}$. Consider the
Hilbert scheme parametrizing all the  hypersurfaces of degree $n$
in $\Pcq$:
\begin{equation}
\label{Hilbert}  \fc\subset \Ps^{2m+1}\times \mid \oc
_{\Ps^{2m+1}}(n) \mid\to \mid \oc _{\Ps^{2m+1}}(n) \mid.
\end{equation}
Notice that by the rational map defined in (\ref{rationalmap}), we
may regard $L$ as a line in $\mid \oc _{\Ps^{2m+1}}(n)\mid$ and
the restriction of the universal family (\ref{Hilbert}) to $L$
gives our pencil $g:\tp\to L$, i.e.  we have:
\begin{equation}
\label{ptilde}  \tp=\fc\times_{\mid \oc _{\Ps^{2m+1}}(n)\mid}L
\subset \Ps^{2m+1}\times \mid \oc _{\Ps^{2m+1}}(n) \mid.
\end{equation}
It follows that the Milnor fibre defined by $g$ in correspondence
of the critical value $a_i\in L$ is equal to the Milnor fibre
defined by  the universal family  (\ref{Hilbert}) at $a_i$.
Therefore, we may interpret the Milnor number $\mu$ of the
singularity of the hypersurface $F_{a_i}$ (recall Remark
(\ref{milnornumber})) as the multiplicity of the discriminant
locus $\dc_n$ of the whole linear system $\mid \oc
_{\Ps^{2m+1}}(n) \mid$ at $a_i$ (see \cite{Loj}, pp. 63-64 and 77,
and \cite{Dimca}, p. 81). Hence, a general line $L' \subset \mid
\oc_{\Ps^{2m+1}}(n) \mid$ passing near to the critical value $a_i$
of $L$ ($a_i\notin L'$) transversally meets $\dc_n$ in  $\mu$
smooth points $b_1,\dots,b_{\mu}$, and the singular locus of each
hypersurface $F_{b_h}$ consists of exactly one  node. Similarly as
in the definition of $I_i$ (see (\ref{defofI})), using this node
we may define a certain cocycle $\delta_{b_h}$ in $H^{2m}(F_t
;\mathbb{C})$. It turns out that these cocycles
$\{\delta_{b_1},\dots,\delta_{b_{\mu}}\}$ lie in $I_i$ and here
they form a distinguished basis (see \cite{Loj}, p. 76, and
\cite{Dimca}, p. 83). And so to prove our Lemma \ref{trivial} it
suffices to prove that the monodromy induced by $L$ trivially acts
on each $\delta_{b_h}$.

To this purpose, fix a $\delta_{b_h}$,  and choose  general germs
$\{F_{\tau}\}_{\tau\in\Delta}$ and $\{G_{\tau}\}_{\tau\in\Delta}$
in $\mid \oc_{\Ps^{2m+1}}(n) \mid$ and in $\mid
\oc_{\Ps^{2m+1}}(k) \mid$, with $F_0=F_t$, $F_{\epsilon}=F_{b_h}$,
$G_0=G$, such that, for any $\tau\neq 0$, $F_{\tau}$ only has one
node belonging also to $G_{\tau}$ (here $\Delta\subset \mathbb{C}$
denotes a closed disk centered at $0$ with small radius
$\epsilon$). Put $X_{\tau}=G_{\tau}\cap F_{\tau}$. Then $X_{\tau}$
is a  complete intersection of dimension $2m-1$ with only one
node, and, as for $X=X_0$, in correspondence of each $X_{\tau}$ we
may define a general pencil $g_{\tau}: \tp_{\tau} \to L_{\tau}$,
with $L_{\tau}\subset \,\mid \ic _{X_{\tau},\Ps^{2m+1}}(n)
\mid\,\subset \,\mid \oc _{\Ps^{2m+1}}(n) \mid$, in such a way
that the family $\{L_{\tau}\}_{\tau\in\Delta}$ trivially deforms
our starting pencil $L=L_0$. Similarly as in the definition of
$I_i$, the Milnor fibre of $g_{\epsilon}$ corresponding to the
nodal fibre $F_{\epsilon}=F_{b_h}=g_{\epsilon}^{-1}(b_h)$, defines
a subspace $I_{\epsilon}\subset
H_{2m}(g_{\epsilon}^{-1}(b_h+\rho_{\epsilon});\mathbb{C})$
($0<\rho_{\epsilon}<<1$). Notice that, via the total space of the
family $\{L_{\tau}\}_{\tau\in\Delta}$, we may transport
$I_{\epsilon}$ in $H_{2m}(F_t ;\mathbb{C})\simeq H^{2m}(F_t
;\mathbb{C})$, and we may assume that in this way we obtain
exactly $Span(\delta_{b_{h}})$. Moreover, using again the total
space of the deformation $\{L_{\tau}\}_{\tau\in\Delta}$, we see
that any closed path in $L-\dc$ is free homotopic in $\mid \oc
_{\Ps^{2m+1}}(n) \mid-\dc_n$ to some closed path contained in
$L_{\epsilon}-\dc_{\epsilon}$ ($\dc_{\epsilon}=$ discriminant
locus of $\mid \ic _{X_{\tau},\Ps^{2m+1}}(n) \mid$). It follows
that if the monodromy of $L_{\epsilon}$ trivially acts on
$I_{\epsilon}$, then also the monodromy of $L$ trivially acts on
$\delta_{b_{h}}$. Therefore, in order to prove Lemma
\ref{trivial}, we may assume $X$ with a unique ordinary double
point.

With this assumption, then $Sing(\qc)=\{q_1,q_{\infty}\}$, and we
only have to prove that the monodromy defined by $L$ trivially
acts on $I_1$. To this purpose, let $\pi\subset \Ps^*$ be a
general projective plane, so that $\pi\cap\hc_1$ is a general line
in $\hc_1$. Denote by $Y$ the set of points in $\pi\cap\hc_1$
parametrizing hyperplanes which are limit of some sequence $z_n$
of tangent hyperplanes at smooth part of $\qc$, such that there
exists a sequence of regular contact points $p_n\in Sing(\qc\cap
H_{z_n})$ converging to $q_1$. Notice that $Y$ is contained in the
finite set $\pi\cap \hc_1\cap\qc^*$. For any $y\in Y$ denote by
$\Delta_{y}$ a closed disk of $\pi\cap \hc_1$, with center $y$ and
positive radius $<<1$, and put $K=(\pi\cap \hc_1)- \bigcup_{y\in
Y}\Delta_y^{\circ}$. Notice also that we may assume our pencil $L$
contained in $\pi$ and close to $\pi\cap \hc_1$, because  such a
pencil is sufficiently general to apply Zariski Theorem, which
ensures that $\pi_1(L-\dc,t)$ maps onto $\pi_1(\Ps^*-\dc,t)$ (see
\cite{La}, (7.4.1), or \cite{Voisin}, Th\'eor\`eme 15.22). Now
consider the restriction of the universal family (\ref{Hilbert})
to $\pi$:
$$
\varphi: \fc_{\pi}:=\fc\times_{\mid \oc
_{\Ps^{2m+1}}(n)\mid}\pi\to \pi.
$$
Recall from (\ref{ptilde}) that we may regard $\tp\subset
\fc_{\pi}$, and the pencil $g:\tp\to L$ as the restriction of
$\varphi$ to $L$. Using \cite{Loj}, Theorem (2.8), we see that for
any $x\in K$ there exists a closed ball $D_{q_1,x}\subset
\Ps^{2m+1}\times \mid \oc _{\Ps^{2m+1}}(n)\mid$, with positive
radius and centered at $q_1$, and a closed ball $C_x \subset \pi$
(with positive radius and of real dimension $4$) centered at $x$,
such that the induced map
$$
\varphi_x: \varphi^{-1}(C_x)\cap D_{q_1,x}\to C_x
$$
is a Milnor fibration whose discriminant locus simply is $C_x \cap
\hc_1$ (observe that
$\pi\cap\hc_1\cap(\qc^*\cup\hc_{\infty})=Y\cup\{pt\}$, with $\mid
Y\mid =2$ and $pt\notin Y$). In view of (\ref{ptilde}), if the
critical value $a_1$ of $L$ corresponding to $q_1$ belongs to some
$C_x$, then the Milnor fibre as defined in (\ref{Milnorfiber})
(note that now we have $i=1$) is  the Milnor fibre of $\varphi_x$.

Moreover, since $x\in K$ then we may assume that for any $z\in C_x
\cap \hc_1$ the map $\varphi_x$ represents the Milnor fibration of
the isolated complete intersection singularity
$(\varphi^{-1}(z)\cap D_{q_1,x}, q_1)$. Therefore, since $K$ is
compact, using the local data $x\in K$, $C_x$ and $D_{q_1,x}$, one
may construct a connected open tubular neighborhood $\mc$ of $K$
in $\pi$, with $a_1\in\mc$, and a closed ball $D_{q_1}\subset
\Ps^{2m+1}\times \mid \oc _{\Ps^{2m+1}}(n)\mid$ of positive radius
and centered at $q_1$ such that the map
\begin{equation}
\label{varphi}  \varphi_{\mc}: z\in\varphi^{-1}(\mc)\cap D_{q_1}
\to\varphi(z)\in \mc
\end{equation}
defines a  $C^{\infty}$-fibre bundle on $\mc -\hc_1$, and whose
fibre $\varphi_{\mc}^{-1}(z)$, $z\in \mc -\hc_1$, may be
identified   with the Milnor fibre of $g:\tp\to L$ corresponding
to $q_1$.

Since also $\pi\cap\hc_1$ is compact,  one may construct an open
tubular neighborhood $\nc$ of $\pi\cap\hc_1$ in $\pi$ in such a
way that one may obtain $\mc$ removing from $\nc$ suitable compact
tubular neighborhoods $\nc_y$ of the disks $\Delta_{y}$, $y\in Y$.

Now, as in \cite{La}, p. 35, Fig. 1, for any critical value $a_j$
of $L$ fix a closed disk $\Delta_j\subset L-\{t\} \simeq \bC$ with
center $a_j$ and  radius $0<\rho <<1$, and, in
$L-\bigcup_{j}\Delta_j^{\circ}$ ($\Delta_j^{\circ}=$ interior of
$\Delta_j$), choose a $C^{\infty}$ path $l_j$ from $t$ to
$a_j+\rho$ with no self-intersection points and such that $l_j\cap
l_h=\{t\}$ for $j\neq h$. Let $w_j\in\pi_1(L-\dc, t)$ be the
homotopy class defined by the path $l_j^{-1}\cdot
\partial\Delta_j \cdot l_j$. These classes span $\pi_1(L-\dc, t)$
with exactly one relation
$$
w_1\cdot w_2\cdots w_s\cdot w_{\infty}=1.
$$
Therefore, to prove Lemma \ref{trivial} it suffices to prove that
any $w_j$, $j\neq 1$, trivially  acts on $I_1$.

From the definition of $\nc$ we may assume that all critical
values of $L$ are contained in $\nc$. Now, with the exception of
the point $\{a_1\}=L\cap \hc_1\cap \mc$, in a neighborhood of  any
critical value $a_j$ of $L$ lying in $\mc$,  the fibration
$\varphi_{\mc}$ is trivial. Since $g:\tp\to L$ is a subbundle of
$\varphi$, and the fibre $\varphi_{\mc}^{-1}(z)$, $z\in \mc
-\hc_1$, identifies with the Milnor fibre of $g$ corresponding to
$q_1$, then $w_j$ trivially acts on $I_1$. Therefore, to complete
the proof of Lemma \ref{trivial}, it remains to analyze the
monodromy corresponding to the set $\dc(L,y)$ of the critical
values of $L$ contained in each $\nc_y$.

To this aim, denote by $\Sigma$ the set of points in $ \hc_1$
which are limit of some sequence $z_n$  of tangent hyperplanes at
smooth part of $\qc$, such that there exists a sequence of regular
contact points $p_n\in Sing(\qc\cap H_{z_n})$ converging to $q_1$.
Since $q_1$ is an ordinary double point for $\qc$, then $\Sigma$
is the dual variety of the tangent cone of $\qc$ at $q_1$. Hence
$\Sigma$ is a nondegenerate irreducible quadric in $\hc_1$. Note
that $Y=\pi\cap \Sigma$, and $|Y|= 2$. Now recall that $X$ is a
complete intersection of a smooth hypersurface in $\Ps^{2m+1}$
with a general hypersurface with a unique ordinary double point.
Hence, for a general point $y\in \Sigma$, $H_y\cap \qc$ has an
isolated singular point at $q_1$ with Milnor number $2$. Combining
(\cite{Loj}, pp. 62-63) with (\cite{Loj}, (5.11.a), p. 77), we
know that for a general line $L'\subset \Ps^*$ passing through
such a point $y$, the intersection multiplicity $m_y(L',\dc)$ of
$L'$ with $\dc$ at $y$ is the sum of the Milnor number of $H_y\cap
\qc$ at $q_1$, with the Milnor number of $\qc$ at $q_1$. Therefore
we have $m_y(L',\dc)=3$. On the other hand, the same argument
applied to the general point $x$ of $\hc_1$, which parametrizes a
hyperplane section with an ordinary double point at $q_1$,  shows
that for the general line $L''\subset \Ps^*$ passing through $x$
one has $m_x(L'',\dc)=2$. It follows that the critical locus of
our pencil $L$ meets $\nc_y$ in exactly one point, corresponding
to some tangent hyperplane to the smooth part of $\qc$.

Since for any $y\in Y$ we have  $|\dc(L,y)|=1$, then our argument
above involving the fibration (\ref{varphi}) proves that $I_1$ is
at least globally invariant under the monodromy action induced by
$L$. Now denote by $\tc(L)$ the set of critical values of $L$
corresponding to tangent hyperplane sections of $\qc$, and fix a
point $a_{j_{0}}\in \tc(L)\cap \mc$. For such a critical value
$a_{j_{0}}$ we just proved that the homotopy class $w_{j_0}$
trivially acts on $I_1$. Therefore, if we denote by
$\delta_{j_{0}}\in H^{2m}(F_t;\mathbb{C})$ the \lq\lq
classical\rq\rq vanishing cocycle generating $V_{j_0}$, from
Picard-Lefschetz formula it follows that for any $\xi\in I_1$ one
has $<\xi,\, \delta_{j_{0}}>=0$. On the other hand, from
(\cite{Voisin}, Proposition 15.23) and (\cite{Loj}, p. 113, Lemma
(7.2)) we know that $\pi_1(L-\qc^*, t)$ irreducibly acts  on the
cocycles determined by tangent hyperplane sections of $\qc$. A
fortiori this holds true for $\pi_1(L-\dc, t)$ and so, from the
global invariance of $I_1$, we deduce that for any $\xi\in I_1$
and any $a_{j}\in \tc(L)$ one has $<\xi,\, \delta_j>=0$. From the
Picard-Lefschetz formula again it follows that for any $a_{j}\in
\tc(L)$ the homotopy class $w_j$ trivially acts on $I_1$. Since we
proved that for any $y\in Y$ one has $\dc(L,y)\subset \tc(L)$,
this concludes the proof of Lemma \ref{trivial}.
\end{proof}
\smallskip

From Lemma \ref{trivial}, (\ref{dec}) and (\ref{VI}) we get
 $V_i\subset I\cap V=0$, and this proves (\ref{van}). In other words we
have:
$$
V=V_{r+1}+\dots+V_s.
$$
This means that $V$ is generated by the vanishing cocycles
determined by the hyperplanes of the pencil $L$ which are tangent
to the smooth part of $Q$. Therefore, as before (see
(\cite{Voisin}, Proposition 15.23) and (\cite{Loj}, p. 113, Lemma
(7.2))) $\pi_1(L-\dc, t)$ irreducibly acts on $V$.
\medskip
\begin{remark}
\label{notasaito} This concludes the proof of Theorem \ref{Saito}.
\end{remark}
\medskip

This enables us to prove that:
\begin{equation}
\label{neron} NS_m(F_t)\subset I
\end{equation}
(as before, we identify $H_{2m} (F_t ; \mathbb{C})\simeq H^{2m}
(F_t ; \mathbb{C})$ via Poincar\`e duality). In fact, argue by
contradiction. Suppose there exists $\xi \in NS_m(F_t)$ such that
${\xi}^{w}\neq \xi$ for some $w\in \pi_1(L-\dc, t)$. We may write
$\xi=i+v$ for some $i\in I$ and $v\in V$, and we have
\begin{equation}
\label{contra} v-v^{w}=\xi - {\xi}^{w}\neq 0.
\end{equation}
Since $\pi_1(L-\dc, t)$ irreducibly acts on $V$, and $NS_m(F_t)$
is globally invariant, (\ref{contra}) implies that $V\subset
H^{m,m}(F_t,\bC)$. On the other hand $\rc$ is birational to
$\Ps^{2m+1}$ and so $H^{2m,0}(\rc,\bC)=0$ (see \cite{Hartshorne},
p. 190, Ex. 8.8). Therefore, since $I=j^*H^4 (\rc ; \mathbb{C})$,
we get
$$
H^{2m,0} (F_t ; \mathbb{C})=I^{2m,0}\oplus V^{2m,0}=0.
$$
This is in contrast with our hypothesis $n>2m+1$. This proves
(\ref{neron}).

We are in position to prove that the natural map $NS_m(X) \to
NS_m(F_t)$ is surjective. To this purpose fix an algebraic class
$\xi \in NS_m(F_t)$, which we may assume represented by some
projective algebraic subvariety $S_1\subset F_t$ of dimension $m$,
and consider the flag Hilbert scheme $\Sc$, with reduced
structure, parametrizing pairs $(S,F)$, with $F\in\lc$ and
$S\subset F$ a projective subvariety of dimension $m$. Let
$\cc\subset \Sc$ be an irreducible projective curve passing
through the point $(S_{1},F_t)$. Since $F_t$ is Noether-Lefschetz
general, we may assume $\cc$ dominating $L$ and such that $t$ is a
regular value of the natural branched covering map $\pi:\cc \to
L$. This curve determines a projective subvariety $T\subset \tp$
of dimension $m+1$, whose intersection with $F_t$ is the union of
all the subvarieties $S_i$, $i=1,\dots,d$, corresponding to the
fibre of $\pi$ over the point $t\in L$ ($d=$ degree of $\pi$). The
monodromy of $\pi$ is transitive, and so by (\ref{neron}) we
deduce that all the $S_i$ are homologous in $F_t$. In other words,
we have
\begin{equation}
\label{csi} \xi={\frac{1}{d}}\cdot j_1^*(T)  \quad{\text{in}}\quad
NS_m(F_t),
\end{equation}
where $j_1^*:A_{m+1}(\tp)\to NS_m(F_t)$ is the natural map induced
by the inclusion $F_t \subset \tp$. Now recall from
(\ref{baselocus}) that $\tp$ is the blowing-up of $\pc$ along the
base locus $B$, which is isomorphic to a projective smooth
complete intersection in $\Ps ^{2m+1}$ of dimension $2m-1$. By
(\cite{Fulton}, p.114-115, Proposition $6.7$, (e)), we know that
$A_{m+1}(B\times {\Ps ^1})\oplus A_{m+1}(\pc)$ maps onto
$A_{m+1}(\tp)$, and therefore we may write
$$
T=lH_{\tp}^m+\alpha^*(Z)\quad {\text{in}}\quad A_{m+1}(\tp),
$$
where $H_{\tp}$ is the pull-back in $\tp$ of the hyperplane class
in ${\Ps ^{2m+1}}$, $l$ is a suitable integer, $\alpha: \tp\to
\pc$ is the natural projection, and $Z$ is a suitable  class in
$A_{m+1}(\pc)$. Plugging previous formula into (\ref{csi}), and
using the natural map $j_2^*:A_{m+1}(\pc)\to NS_m(F_t)$ induced by
the inclusion $F_t \subset \pc$, we get
\begin{equation}
\label{geidue} j_1^*(T) = j_2^*(lH_{\pc}^m+Z),
\end{equation}
where $H_{\pc}$ is the pull-back in $\pc$ of the hyperplane class
in ${\Ps ^{2m+1}}$. Since $\pc$ is the blowing-up of $\Ps ^{2m+1}$
along $X$, again by (\cite{Fulton}, l.c.) we know that
$A_{m+1}(\tilde X)\oplus A_{m+1}(\Ps ^{2m+1})$ maps onto
$A_{m+1}(\pc)$, where $h:\tilde X\subset \pc$ is the exceptional
divisor, which in turn is a $\Ps ^1-$bundle $\beta: \tilde X\to X$
over the complete intersection $X$. The group $A_{m+1}(\tilde X)$
is spanned by $\beta^*(A_m(X))$, and by the cycles obtained
intersecting a fixed section of $\beta$ with
$\beta^*(A_{m+1}(X))$. As a  section we may choose $\tilde G\cap
\tilde X$ (compare with (\ref{gtilde})). It follows that we may
write
$$
lH_{\pc}^m+Z= lH_{\pc}^m +h_{*}\beta^*(W_1)+h_{*}(\tilde G\cap
\tilde X\cap \beta^{*}(W_2)) \quad {\text{in}}\quad A_{m+1}(\pc),
$$
where $W_1$ and $W_2$ are suitable classes in $A_m(X)$ and
$A_{m+1}(X)$. Taking into account that $\tilde G\cap \tilde X$ is
disjoint with $F_t$, from (\ref{csi}) and (\ref{geidue}) we obtain
$$
\xi={\frac{1}{d}}\cdot j_2^*(lH_{\pc}^m
+h_{*}\beta^*(W_1)+h_{*}(\tilde G\cap \tilde X\cap
\beta^{*}(W_2)))={\frac{1}{d}}\cdot
\gamma_{*}({\frac{l}{k}}H^{m-1}_X+W_1),
$$
where $H_X$ is the hyperplane section of $X$, $k=deg(G)$ and
$\gamma_{*}$ denotes the  map $NS_m(X) \to NS_m(F_t)$. This proves
that this map is onto, and concludes the proof of Theorem
\ref{mmain}.
\end{proof}

\bigskip
{\it Proof of Theorem \ref{NL}}. We are in position to prove
Theorem \ref{NL} stated in the Introduction. To this purpose
consider a general pencil of hypersurface sections $\{\qc\cap
R_t\}_{t\in\Ps^1}$ of $\qc$, with $deg (R_t)>>0$. With the same
methods we used in the proof of Theorem \ref{Saito}, for a general
$t$, we may prove an orthogonal decomposition $H^{2m} (\qc\cap
R_t; \mathbb{C})=I\oplus V$, such that the monodromy
representation of the pencil is irreducible on $V$, and
$I=j^*H^{2m}(\rc; \mathbb{C})$, where $\rc$ denotes  a certain
desingularization of $\qc$ and $j$  the inclusion $\qc\cap
R_t\subset\rc$. Since $deg (R_t)>>0$ then $h^{2m,0} (\qc\cap R_t)>
h^{2m,0}(\rc)$ (racall that $h^{2m,0}(\rc)$ is a birational
invariant, and so it only depends on $\qc$). It follows that
$NS_m(\qc\cap R_t)\subset I$. A similar argument as in the proof
of (\ref{csi}) shows that $NS_{m+1}(\rc)$ maps onto $NS_m(\qc\cap
R_t)$. On the other hand, since $R_t$ does not meet the singular
locus of $\qc$, then for the Gysin morphisms $a: H_{2m+2}(\rc;
\mathbb{C})\to H_{2m}(\qc\cap R_t; \mathbb{C})$ and $b:
H_{2m+2}(\qc; \mathbb{C})\to H_{2m}(\qc\cap R_t; \mathbb{C})$ one
has $a=b\circ p$, where $p: H_{2m+2}(\rc; \mathbb{C})\to
H_{2m+2}(\qc; \mathbb{C})$ denotes the push-forward. Therefore
also $NS_{m+1}(\qc)$ maps onto $NS_m(\qc\cap R_t)$, and so
$dim(NS_m(\qc\cap R_t))=1$. This concludes the proof of Theorem
\ref{NL}.
$\qquad\qquad\qquad\qquad\qquad\qquad\qquad\qquad\qquad\qquad
\square$
\medskip

\medskip
{\it Proof of Theorem \ref{Saitodue}}. With the same notation as
in the proof of Theorem \ref{mmain},  we are going to prove that
the image $I_X$ of $H_{2m}(X; \mathbb{C})$ in $H_{2m} (F;
\mathbb{C})$ $\simeq$ $H^{2m}(F; \mathbb{C})$ is equal to $I$.
First notice that $I_X\subset I$ because the cycles coming from
$X$ are invariant. So it suffices to prove that $I\subset I_X$.
Since $I=j^*H^{2m}(\rc ; \mathbb{C})$, via Poincar\`e duality we
see that $I$ is equal to the image of the Gysin morphism $a:
H_{2m+2}(\rc ; \mathbb{C})\to H_{2m}(F; \mathbb{C})$. Since $F$,
as subvariety of $\pc$, does not meet the singular locus of $\pc$,
then one has $a=b\circ p$, where $p: H_{2m+2}(\rc; \mathbb{C})\to
H_{2m+2}(\pc; \mathbb{C})$ denotes the push-forward and $b$ the
Gysin morphism $H_{2m+2}(\pc; \mathbb{C})\to H_{2m}(F;
\mathbb{C})$. Therefore $I$ is contained in the image of
$H_{2m+2}(\pc; \mathbb{C})$ through $b$. Now denote by $\tilde X$
the exceptional divisor of $\pc$. From \cite{La}, p. 23, we know
there  exists a natural isomorphism $H_{*} (\pc,\tilde X;
\mathbb{C})$ $\simeq$  $H_{*}(\Ps^{2m+1}, X; \mathbb{C})$. On the
other hand, using \cite{Dimca}, Theorem 4.3, p. 161, one sees that
$H_{2m+2}(\Ps^{2m+1}, X; \mathbb{C})=H_{2m+3}(\Ps^{2m+1}, X;
\mathbb{C})=0$. Hence the inclusion $\tilde X\subset \pc$ induces
a natural isomorphism $H_{2m+2}(\pc; \mathbb{C})\simeq
H_{2m+2}(\tilde X; \mathbb{C})$, and so $H_{2m+2}(\tilde X;
\mathbb{C})$ maps onto $I$. Taking into account that $\tilde X$ is
a $\Ps^1$-bundle over $X$, from Leray-Hirsch Theorem
(\cite{Spanier}, p.  258) we know that all the homology of $\tilde
X$ comes from $X$, up the cycles contained in a fixed section of
the bundle $\tilde X\to X$, which we may choose disjoint with $F$.
Therefore $I_X$ contains $I$, and  this concludes the proof of
Theorem \ref{Saitodue}.
$\qquad\qquad\qquad\qquad\qquad\qquad\qquad\qquad\qquad\qquad\qquad
\qquad\qquad\qquad\qquad\qquad\quad \square$
\medskip

\bigskip
\bigskip

{\bf{Aknowledgements}}
\smallskip

We would like to thank Ciro Ciliberto for valuable discussions and
suggestions on the subject of this paper.


\bigskip
\bigskip


\begin{thebibliography}{subsurfPq}



\bibitem{BS} Beltrametti, M.C. - Sommese, A.J.: {\it The adjunction theory of complex projective
varieties}, de Gruyter Expositions in Mathematics, 1995.

\bibitem{Cheltsov1} Cheltsov, I.: {\it On factoriality of nodal
threefolds}, J. Alg. Geom. {\bf 14}, n. 4, 663-690 (2005).

\bibitem{Cheltsov2} Cheltsov, I.: {\it Factorial nodal threefolds in
$\Pcq$}, AG/0410252, (2004).

\bibitem{Cheltsov3} Cheltsov, I.: {\it Points in projective space and
applicatipons}, AG/0511578, (2005).

\bibitem{Cheltsov4} Cheltsov, I. - Park, J.: {\it Factorial hypersurfaces in $\mathbb{P}^4$ with nodes},
AG/0511673, (2005).

\bibitem{CD1} Ciliberto, C. - Di Gennaro V.: {\it Factoriality of certain
hypersurfaces of $\mathbb{P}^4$ with ordinary double points},
Encyclopaedia of Mathematical Sciences, Springer-Verlag, Vol. 132,
1-7 (2004).

\bibitem{CD2} Ciliberto, C. - Di Gennaro V.: {\it Factoriality of certain threefolds complete intersection in
$\mathbb{P}^5$ with ordinary double points}, Comm. in Alg. {\bf
32}, No. 7,  2705-2710 (2004).


\bibitem{Clemens} Clemens, C. H.: {\it Double solids}, Adv. in Math.
{\bf 47}, 107-230, (1983).


\bibitem{Cynk} Cynk, S.: {\it Defect of a nodal hypersurface},
Man. Math. {\bf 104}, 325-331 (2001).

\bibitem{Dimca} Dimca, A.: {\it Singularity and Topology of
Hypersurfaces}, Springer Universitext, New York 1992.



\bibitem{Vogel} Flenner, H. - O'Carroll L. - Vogel W.: {\it Joins and
intersections}, Springer-Verlag, Berlin, 1999.


\bibitem{Fulton} Fulton, W.: {\it Intersection theory}, Ergebnisse
der Mathematik und ihrer Grenzgebiete; 3.Folge, Bd. 2,
Springer-Verlag 1984.

\bibitem{GM} Geramita, A.V. - Maroscia, P.: {\it The ideal of forms
vanishing at a finite set of points in $\Ps^n$}, J. Algebra {\bf
90}, 528-555 (1984).



\bibitem{Greco} Greco, S.: {\it Normal varieties},
Published by Istituto Nazionale di Alta matematica Roma,
Distributed by Academic Press London and New York, 1978.

\bibitem{Hartshorne} Hartshorne, R.: {\it Algebraic Geometry},
Graduate Texts in Mathematics, 52, Springer-Verlag, 1983.


\bibitem{Hu} Hu, W.: {\it The Generalized Hodge conjecture for 1-cycles and codimension two
algebraic cycles}, AG/0511725, (2005).


\bibitem{KM} Koll\'ar, J. - Mori, S.: {\it Birational Geometry of Algebraic
Varieties}, Cambridge Tracts in Mathematics 134, Cambridge
University Press, Cambridge, 1998.

\bibitem{La} Lamotke, K.: {\it The topology of complex projective varieties
after S. Lefschetz}, Topology {\bf 20}, 15-51 (1981).


\bibitem{Lazarsfeld} Lazarsfeld, R.: {\it Positivity in algebraic geometry.
II.}, Ergebnisse der Mathematik und ihrer Grenzgebiete, 3.Folge. A
Series of Modern Surveys in Mathematics, 49, Springer-Verlag,
Berlin, 2004.



\bibitem{Loj} Looijenga, E.J.N.: {\it Isolated Singular Points on Complete
Intersections}, London Mathematical Society Lecture Note Series
{\bf 77}, Cambridge University Press, Cambridge, 1984.

\bibitem{Lopez} Lopez, A. F.: {\it Noether-Lefschetz theory and the Picard group of projective surfaces},
Mem. Amer. Math. Soc.  {\bf 89}, (1991).

\bibitem{Ml} Mella, M.: {\it Birational geometry of quartic
3-folds II: the importance of being $\mathbb{Q}$-factorial}, Math.
Ann. $\bf{330}$, 107-126 (2004).

\bibitem{MP} Miyaoka, Y. - Peternell, T.: {\it Geometry of Higher Dimensional Algebraic
Varieties}, DMV Seminar, Bd. {\bf 26} (Birkh\"auser, Basel 1997).

\bibitem{OS} Otwinowska, A. - Saito, M.: {\it Monodromy of a family of hypersurfaces containing a given subvariety},
Ann. Scient. \'Ec. Norm. Sup., $4^{e}$ s\'erie, t. 38, 365-386
(2005).

\bibitem{PS} Parshin, A.N. - Shafarevich, I.R. (Eds.): {\it Algebraic Geometry III},
Enciclopaedia of Math. Sciences, vol. 36, Springer, Berlin, 1998.



\bibitem{Spanier} Spanier, E.H.: {\it Algebraic Topology}, McGraw-Hill
Series in Higher Mathematics, 1966

\bibitem{Voisin} Voisin, C.: {\it Th\'eorie de Hodge et g\'eom\'etrie alg\'ebrique
complexe}, Cours Sp\'ecialis\'es 10, Soci\'et\'e Math\'ematique de
France 2002.

\bibitem{Werner} Werner, J.: {\it Kleine Aufl\"osungen spezieller dreidimensionaler
Variet\"aten}, Bonner Mathematische Schriften {\bf 186} (1987),
Universit\"at Bonn, Mathematisches Institut, Bonn.


\end{thebibliography}
\end{document}